\documentclass{article}

\usepackage{amsmath,amscd,amsfonts,amssymb,amsthm}
\usepackage{subfig}
\usepackage{graphicx}
\usepackage{epstopdf}
\usepackage{a4wide}

\newcommand{\subj}[1]{\par\noindent{\bf Mathematics Subject Classification 2010: }#1.}
\newcommand{\keyw}[1]{\par\noindent{\bf Keywords: }#1.}

\theoremstyle{definition}
\newtheorem{definition}{Definition}

\theoremstyle{remark}

\def\a{\alpha}
\def\LDa{{_aD_t^\a}}
\def\LCD{{_a^CD_t^\a}}
\def\DS{\displaystyle}

\begin{document}

\title{A Discrete Method to Solve Fractional\\ Optimal Control Problems\thanks{This is 
a preprint of a paper whose final and definite form will appear in 
\emph{Nonlinear Dynamics}, ISSN 0924-090X. Paper submitted 26/Nov/2013; 
revised 17/March/2014; accepted for publication 20/March/2014.}}

\author{Ricardo Almeida\\
{\tt ricardo.almeida@ua.pt}
\and
Delfim F. M. Torres\\
{\tt delfim@ua.pt}}

\date{Center for Research and Development in Mathematics and Applications (CIDMA)\\
Department of Mathematics, University of Aveiro, 3810--193 Aveiro, Portugal}

\maketitle


\begin{abstract}
We present a method to solve fractional optimal control problems,
where the dynamic depends on integer and Caputo fractional derivatives.
Our approach consists to approximate the initial fractional order problem
with a new one that involves integer derivatives only.
The latter problem is then discretized,
by application of finite differences, and solved numerically.
We illustrate the effectiveness of the procedure with an example.
\end{abstract}

\subj{26A33, 49M25}

\keyw{fractional calculus, fractional optimal control, direct methods}


\section{Introduction}

The fractional calculus deals with differentiation
and integration of arbitrary (noninteger) order \cite{book:samko}.
It has received an increasing interest due to the fact that fractional
operators are defined by natural phenomena \cite{Benson}.

A particularly interesting, and very active, research area is that of
the fractional optimal control, where the dynamic control system involves
not only integer order derivatives but fractional operators as well \cite{book:frac,pooseh:press}.
Such fractional problems are difficult to solve and numerical methods are often applied,
which have proven to be efficient and reliable \cite{sug:r3,ford1,ford2,Lotfi,Almeida1}.
Typically, using an approximation for the fractional operators, the fractional optimal control
problem is rewritten into a new one, that depends on integer derivatives only \cite{Atanackovic,Almeida}.
Then, using necessary optimality conditions, the problem is reduced
to a system of ordinary differential equations and, by finding its solution,
one approximates the solution to the original fractional problem \cite{Almeida}.
Our approach is different here: after replacing the fractional operator,
we consider the augmented functional and apply an Euler-like method.

The text is organized as follows. In Section~\ref{sec:2}
we briefly recall the necessary concepts and results.
Our method is presented in Section~\ref{sec:3}
and then illustrated, with an example, in Section~\ref{sec:4}.


\section{Preliminaries}
\label{sec:2}

We start with some definitions needed in the sequel.
For more on the subject we refer the reader
to \cite{Kilbas,book:frac,Podlubny,book:samko}.

\begin{definition}
\label{def:DF:IF}
Let $x:[a,b]\to\mathbb R$ and $\alpha>0$
be the order of the integral or the derivative.
For $t \in[a,b]$, we define
\begin{enumerate}
\item the left Riemann--Liouville fractional integral by
$$
{_aI_t^{\alpha}}x(t)=\frac{1}{\Gamma(\alpha)}
\int_a^t (t-\tau)^{\alpha-1}x(\tau) d\tau,
$$
\item the right Riemann--Liouville fractional integral by
$$
{_tI_b^{\alpha}}x(t)=\frac{1}{\Gamma(\alpha)}
\int_t^b (\tau-t)^{\alpha-1}x(\tau) d\tau,
$$
\item the left Riemann--Liouville fractional derivative by
$$
{_aD_t^{\alpha}}x(t)=\frac{1}{\Gamma(n-\alpha)}
\frac{d^n}{dt^n}\int_a^t (t-\tau)^{n-\alpha-1}x(\tau) d\tau,
$$
\item the right Riemann--Liouville fractional derivative by
$$
{_tD_b^{\alpha}}x(t)=\frac{(-1)^n}{\Gamma(n-\alpha)}
\frac{d^n}{dt^n}\int_t^b (\tau-t)^{n-\alpha-1}x(\tau) d\tau,
$$
\item the left Caputo fractional derivative by
$$
{_a^CD_t^{\alpha}}x(t)=\frac{1}{\Gamma(n-\alpha)}
\int_a^t (t-\tau)^{n-\alpha-1}x^{(n)}(\tau) d\tau,
$$
\item the right Caputo fractional derivative by
$$
{_t^CD_b^{\alpha}}x(t)=\frac{(-1)^n}{\Gamma(n-\alpha)}
\int_t^b (\tau-t)^{n-\alpha-1}x^{(n)}(\tau) d\tau,
$$
\end{enumerate}
where $n=[\alpha]+1$ for the definitions of
Riemann--Liouville fractional derivatives, and
\begin{equation}
\label{def:n}
n=\left\{\begin{array}{lll}
[\alpha]+1 & \mbox{ if } \, &\a \notin \mathbb N_0,\\
\alpha & \mbox{ if } \, &\a \in \mathbb N_0,\\
\end{array}\right.
\end{equation}
for the Caputo fractional derivatives.
\end{definition}

When $\alpha=n$ is an integer, the operators
of Definition~\ref{def:DF:IF}
reduce to standard ones:
\begin{equation*}
\begin{split}
{_aI_t^{n}}x(t)
&= \displaystyle \int_a^t d\tau_1 \int_a^{\tau_1} d\tau_2
\ldots \int_a^{\tau_{n-1}}x(\tau_n) d\tau_n,\\
{_tI_b^{n}}x(t) &= \displaystyle\int_t^b d\tau_1
\int_{\tau_1}^b d\tau_2 \ldots \int_{\tau_{n-1}}^b x(\tau_n) d\tau_n,\\
{_aD_t^{n}}x(t) &= \displaystyle  {_a^CD_t^{n}}x(t)  \,
=  \, x^{(n)}(t),\\
{_tD_b^{n}}x(t) &= \displaystyle {_t^CD_b^{n}}x(t)
\, = \, (-1)^n x^{(n)}(t).
\end{split}
\end{equation*}

There is an useful relation between Riemann--Liouville and Caputo
fractional derivatives:
\begin{equation}
\label{relation}
{_a^CD_t^\alpha}x(t)={_aD_t^\alpha}x(t)-\sum_{k=0}^{n-1}
\frac{x^{(k)}(a)}{\Gamma(k-\alpha+1)}(t-a)^{k-\alpha}.
\end{equation}

For numerical purposes, approximations are used to deal with these fractional operators.
The Riemann--Liouville derivatives are expandable in a power series involving integer order derivatives only.
If $x$ is an analytic function, then (cf. \cite{Kilbas})
\begin{equation}
\label{app1}
{_aD_t^{\alpha}} x(t)=\sum_{k=0}^\infty \binom{\a}{k}
\frac{(t-a)^{k-\a}}{\Gamma(k+1-\a)}x^{(k)}(t),
\end{equation}
where $$\binom{\a}{k}=\frac{(-1)^{k-1}\a\Gamma(k-\a)}{\Gamma(1-\a)\Gamma(k+1)}.$$
The obvious disadvantage of using \eqref{app1} in numerical computations
is that in order to have a small error, one has to sum a large number
of terms and thus the function has to possess higher-order derivatives,
which is not suitable for optimal control. To address this problem,
a second approach was carried out in \cite{Atanackovic,Almeida},
where a good approximation is obtained without the requirement
of such higher-order smoothness on the admissible functions.
The method can be explained, for left-derivatives, in the following way.
Let $\a\in(0,1)$ and $x\in C^2[a,b]$. Then,
\begin{equation*}
\LDa x(t)=A(\a)(t-a)^{-\a}x(t)+B(\a)(t-a)^{1-\a}\dot{x}(t)
-\sum_{p=2}^{\infty}C(\a,p)(t-a)^{1-p-\a}V_p(t),
\end{equation*}
where
$V_p(t)$ is the solution of the system
$$
\begin{cases}
\dot{V}_p(t)=(1-p)(t-a)^{p-2}x(t),\\
V_p(a)=0,
\end{cases}
$$
for $p=2,3,\ldots,$ and $A$, $B$ and $C$ are given by
\begin{equation*}
\begin{split}
A(\a) &= \frac{1}{\Gamma(1-\a)}
\left[1+\sum_{p=2}^{\infty}\frac{\Gamma(p-1+\a)}{\Gamma(\a)(p-1)!}\right],\\
B(\a) &= \frac{1}{\Gamma(2-\a)}
\left[1+\sum_{p=1}^{\infty}\frac{\Gamma(p-1+\a)}{\Gamma(\a-1)p!}\right],\\
C(\a,p) &= \frac{1}{\Gamma(2-\a)\Gamma(\a-1)}\frac{\Gamma(p-1+\a)}{(p-1)!}.
\end{split}
\end{equation*}
Using \eqref{relation}, a similar formula can be deduced
for the Caputo fractional derivative. When we consider
finite sums only, that is, when we use the approximation
\begin{equation}
\label{approximation}
\LDa x(t)\approx A(\a,K)(t-a)^{-\a}x(t)+B(\a,K)(t-a)^{1-\a}\dot{x}(t)
-\sum_{p=2}^{K}C(\a,p)(t-a)^{1-p-\a}V_p(t),
\end{equation}
where $K\geq 2$ and
\begin{equation*}
\begin{split}
A(\a,K) &= \frac{1}{\Gamma(1-\a)}\left[1+\sum_{p=2}^{K}
\frac{\Gamma(p-1+\a)}{\Gamma(\a)(p-1)!}\right],\\
B(\a,K) &= \frac{1}{\Gamma(2-\a)}\left[1+\sum_{p=1}^{K}
\frac{\Gamma(p-1+\a)}{\Gamma(\a-1)p!}\right],
\end{split}
\end{equation*}
the error is bounded by
\begin{equation*}
|E_{tr}(t)|\leq \max_{\tau \in [a,t]}\left|\ddot{x}(\tau)\right|
\frac{\exp((1-\a)^2+1-\a)}{\Gamma(2-\a)(1-\a)K^{1-\a}}(t-a)^{2-\a}.
\end{equation*}
See \cite{Almeida} for proofs and other details.


\section{Problem Statement and the Approximation Method}
\label{sec:3}

The fractional optimal control problem that we consider here is the following one.
Let $\a\in(0,1)$ be the fixed fractional order, and consider
two differentiable functions  $L$ and $f$ with domain
$[a,b]\times \mathbb{R}^2$. Minimize the functional
$$
J(x,u)=\int_a^b L(t,x(t),u(t))\,dt
$$
subject to the fractional dynamic constraint
$$
M \dot{x}(t) + N\LCD x(t)
= f\left(t,x(t),u(t)\right), \quad t\in[a,b],
$$
and the initial condition
$$
x(a)=x_a,
$$
where $(M,N)\not=(0,0)$ and $x_a$ is a fixed real number.
Two situations are considered: $x(b)$ fixed or free.
Sufficient and necessary conditions to obtain solutions
for this problem were studied in \cite{pooseh:press}.
Here we proceed with a different approach. First, we replace
the operator $\LCD x(t)$ with the approximation given in \eqref{approximation}.
With relation \eqref{relation} we get
\begin{multline*}
M \dot{x}(t) + N\Bigg[A(t-a)^{-\a}x(t)+B(t-a)^{1-\a}\dot{x}(t)\\
-\sum_{p=2}^{K}C_p(t-a)^{1-p-\a}V_p(t)-\frac{x_a(t-a)^{-\a}}{\Gamma(1-a)}\Bigg]
= f\left(t,x(t),u(t)\right),
\end{multline*}
where, for simplicity,
$$
A=A(\a,K)\, , \quad B= A(\a,K) \, \quad \mbox{and} \quad C_p=C(\a,p).
$$
Thus, one has
$$
\dot{x}(t)
=\frac{f(t,x(t),u(t))-NA(t-a)^{-\a}x(t)
+\DS\sum_{p=2}^K NC_p(t-a)^{1-p-\a}V_p(t)
+\DS\frac{Nx_a(t-a)^{-\a}}{\Gamma(1-a)}}{M+NB(t-a)^{1-\a}}.
$$
Define the vector
$$
\overline V(t)=\left(V_2(t),V_3(t),\ldots,V_K(t)\right)
$$
and the new function
$$
F(t,x,\overline V,u)
=\frac{f(t,x,u)-NA(t-a)^{-\a}x
+\DS\sum_{p=2}^K NC_p(t-a)^{1-p-\a}V_p
+\DS\frac{Nx_a(t-a)^{-\a}}{\Gamma(1-a)}}{M+NB(t-a)^{1-\a}}.
$$
We construct a new optimal control problem: minimize the functional
\begin{equation}
\label{new1}
\overline{J}(x,\overline V,u)=\int_a^b L(t,x(t),u(t))\,dt
\end{equation}
\label{new2}
subject to the dynamic constraints
\begin{equation}
\begin{cases}
\dot{x}(t)=F\left(t,x(t),\overline V(t),u(t)\right),\\[0.25cm]
\dot{V}_p(t)=(1-p)(t-a)^{p-2}x(t), \quad p=2,\ldots,K,
\end{cases}
\end{equation}
and the initial conditions
\begin{equation}
\label{new3}
\begin{cases}
x(a)=x_a,\\
V_p(a)=0, \quad p=2,\ldots,K.
\end{cases}
\end{equation}
To solve the problem \eqref{new1}--\eqref{new3},
one can consider the Hamiltonian function
$$
H(t,x,\overline V,\overline\lambda,u)=L(t,x,u)
+\lambda_1 F(t,x,\overline V,u)+\sum_{p=2}^K \lambda_p (1-p)(t-a)^{p-2}x,
$$
where $\overline\lambda$ denotes the Lagrange multipliers,
$$
\overline\lambda(t)=(\lambda_1(t),\lambda_2(t),\ldots,\lambda_K(t)).
$$
By the Pontryagin maximum principle \cite{Pontryagin},
to solve the problem one should solve the following system of ODEs:
$$
\begin{cases}
\DS\frac{\partial H}{\partial u}=0,\\[0.25cm]
\DS\frac{\partial H}{\partial \lambda_1}=\dot{x},\\[0.25cm]
\DS\frac{\partial H}{\partial \lambda_p}=\dot{V}_p, \quad p=2,\ldots,K,\\[0.25cm]
\DS\frac{\partial H}{\partial x}=-\dot{\lambda}_1,\\[0.25cm]
\DS\frac{\partial H}{\partial V_p}=-\dot{\lambda}_p,\quad p=2,\ldots,K,
\end{cases}
$$
subject to the boundary conditions
$$
\left\{ \begin{array}{l}
x(a)=x_a,\\
V_p(a)=0, \quad p=2,\ldots,K,\\
\lambda_p(b)=0, \quad p=1,\ldots,K,
\end{array}\right.
$$
if $x(b)$ is free, or
$$
\left\{ \begin{array}{l}
x(a)=x_a,\\
V_p(a)=0, \quad p=2,\ldots,K,\\
\end{array}\right.
$$
otherwise. Instead of this indirect approach, we apply here a direct method,
based on an Euler discretization, to obtain a finite-dimensional approximation
of the continuous problem \eqref{new1}--\eqref{new3}.
We summarize briefly the method. Consider a finite grid
$$
a=t_0<t_1<\cdots<t_n=b
$$
where, for simplicity, $t_{i+1}-t_{i}:=\Delta t$ is assumed constant
for all $i\in\{0,\ldots,n-1\}$. Each dynamic constraint is approximated by
$$
x_{i+1}=x_i+\Delta t F(t_i,x_i,\overline V_i,u_i)
$$
and
$$
V_{p,i+1}=V_{p,i}+\Delta t (1-p)(t_i-a)^{p-2}x_i,
\quad p=2,\ldots,K,
$$
where
$$
x_i=x(t_i), \quad u_i=u(t_i), \quad \overline V_i
=(V_2(t_i),V_3(t_i),\ldots,V_K(t_i))
\quad \mbox{and} \quad V_{p,i}=V_p(t_i).
$$
The integral
$$
\int_a^b L\left(t,x(t),u(t)\right)\,dt
$$
in replaced by the Riemann sum
$$
\Delta t \sum_{i=0}^{n-1} L\left(t_i,x_i,u_i\right).
$$
The finite version of problem \eqref{new1}--\eqref{new3}
is the following one: minimize
\begin{equation*}
\Delta t \sum_{i=0}^{n-1} L\left(t_i,x_i,u_i\right)
\end{equation*}
subject to the dynamic constraints
\begin{equation*}
\begin{cases}
x_{i+1}=x_i+\Delta t F\left(t_i,x_i,\overline V_i,u_i\right),\\[0.25cm]
V_{p,i+1}=V_{p,i}+\Delta t (1-p)(t_i-a)^{p-2}x_i, \quad p=2,\ldots,K,
\end{cases}
\end{equation*}
and the initial conditions
\begin{equation*}
\begin{cases}
x_0=x_a,\\
V_{p,0}=0, \quad p=2,\ldots,K
\end{cases}
\end{equation*}
($x(b)$ fixed or free). In the next section we illustrate the method with an example.


\section{An Illustrative Example}
\label{sec:4}

\begin{figure}[!htb]
\centering
\subfloat[\footnotesize{State function $x(t)$.}]{\label{plotX}
\includegraphics[scale=0.30,angle=270]{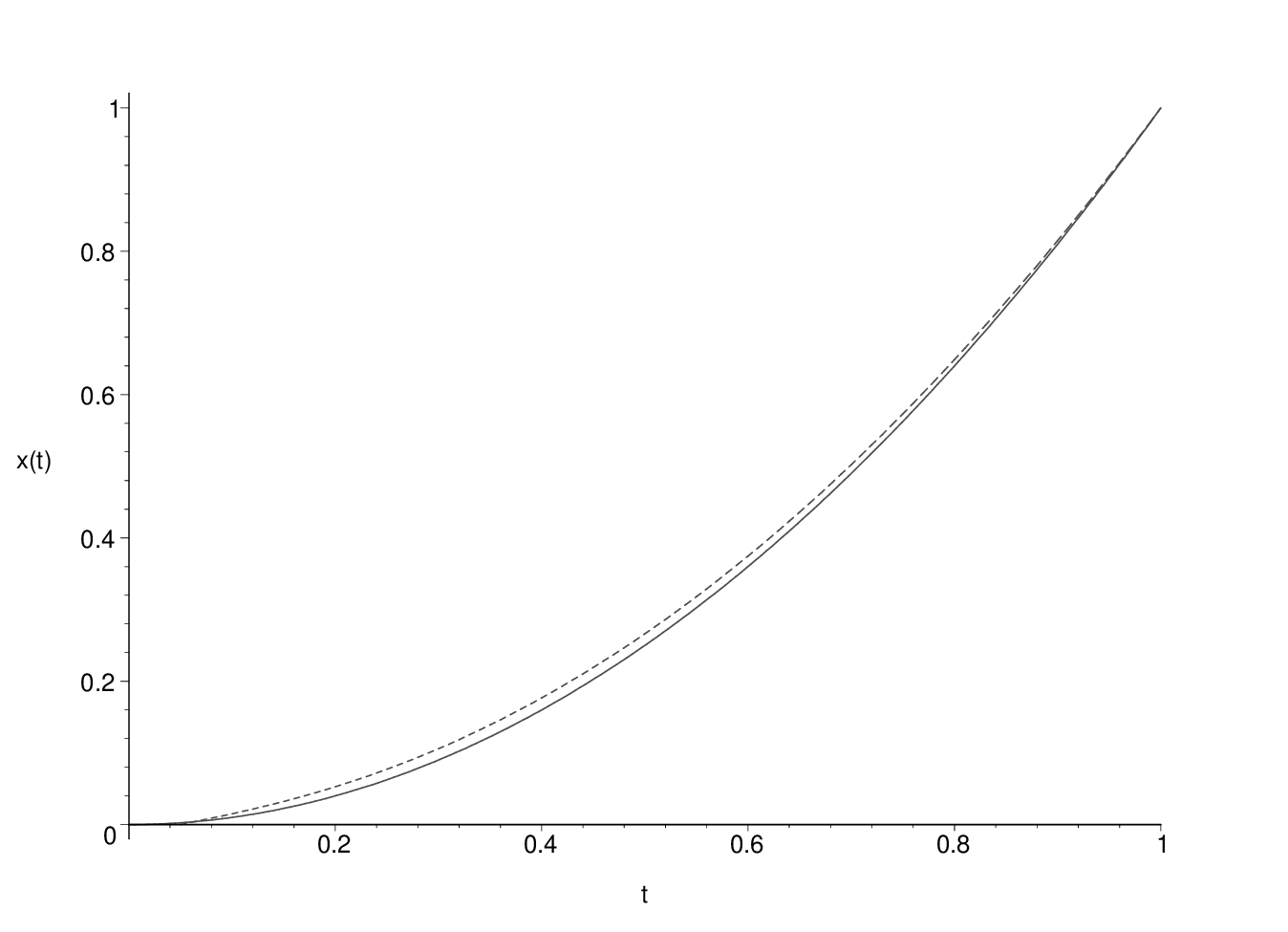}}
\subfloat[\footnotesize{Control function $u(t)$.}]{\label{plotU}
\includegraphics[scale=0.30,angle=270]{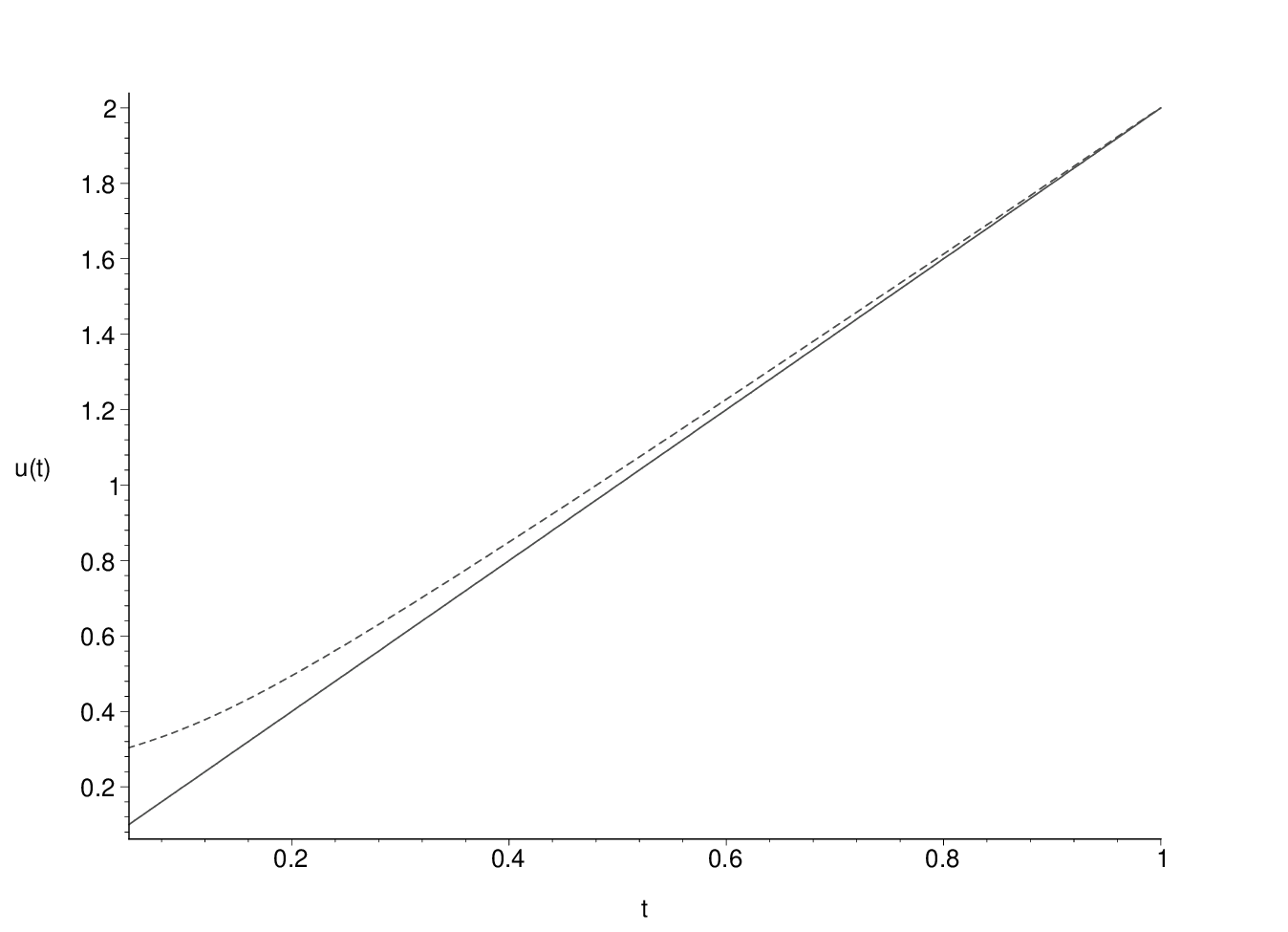}}
\caption{Exact solution to \eqref{ex1}--\eqref{eq:bc:ex} (continuous line)
versus numerical approximation to solution of \eqref{eq:J:ap}--\eqref{eq:bc:ap}
with $K=3$ and $n=100$ (dot line).}
\label{fig:X:U}
\end{figure}

In this section we exemplify the procedure of Section~\ref{sec:3} with a concrete example.
To start, we recall the Caputo fractional derivative of a power function
(cf. Property~2.16 of \cite{Kilbas}). Let $\alpha>0$ be arbitrary, $\beta>n$
with $n$ given by \eqref{def:n}, and define $x(t)=(t-a)^{\beta-1}$. Then,
$$
\LCD x(t) = \frac{\Gamma(\beta)}{\Gamma(\beta-\alpha)}(t-a)^{\beta-\alpha-1}.
$$
Assume that we wish to minimize the functional
\begin{equation}
\label{ex1}
J(x,u)=\int_0^1 (u^2(t)-4x(t))^2\,dt
\end{equation}
subject to the dynamic constraint
\begin{equation}
\label{eq:dc:ex}
\dot{x}(t) +{_0^CD_t^{0.5}} x(t)
=u(t)+\frac{2}{\Gamma(2.5)}t^{1.5}, \quad t\in[0,1],
\end{equation}
and the boundary conditions
\begin{equation}
\label{eq:bc:ex}
x(0)=0 \quad \mbox{and} \quad x(1)=1.
\end{equation}
The solution is the pair
\begin{equation}
\label{eq:sol:ex}
(\overline{x}(t),\overline{u}(t))=(t^2,2t).
\end{equation}
Indeed, \eqref{eq:sol:ex} satisfies both constraints \eqref{eq:dc:ex} and \eqref{eq:bc:ex}
with $J(\overline{x},\overline{u}) = 0$ while functional \eqref{ex1} is non-negative:
$J(x,u) \ge 0$ for all admissible pairs $(x,u)$. After replacing the fractional
derivative by the appropriate approximation, and fixing $K\geq2$,
we get the following problem: minimize the functional
\begin{equation}
\label{eq:J:ap}
\overline{J}(x,\overline V,u)=\int_0^1 (u^2(t)-4x(t))^2\,dt
\end{equation}
subject to the dynamic constraints
\begin{equation}
\label{eq:dc:ap}
\begin{cases}
\dot{x}(t)=\displaystyle \frac{u(t)+\frac{2}{\Gamma(2.5)}t^{1.5}
-At^{-0.5}x(t)+\sum_{p=2}^K C_pt^{0.5-p}V_p(t)}{1+Bt^{0.5}},\\[0.25cm]
\dot{V}_p(t)=(1-p)t^{p-2}x(t), \quad p=2,\ldots,K,
\end{cases}
\end{equation}
and the boundary conditions
\begin{equation}
\label{eq:bc:ap}
\begin{cases}
x(0)=0,\\
V_p(0)=0, \quad p=2,\ldots,K,\\
x(1)=1.\\
\end{cases}
\end{equation}
We apply a direct method to the previous problem, also called a ``first discretize then optimize method''.
More precisely, the objective functional \eqref{eq:J:ap} and the system of ODEs  \eqref{eq:bc:ap}
were discretized by a simple Euler method with fixed step size, as explained in Section~\ref{sec:3}.
Then we solved the resulting nonlinear programming problem by using the AMPL
modeling language for mathematical programming \cite{AMPL}
in connection with IPOPT of W\"achter and Biegler \cite{IPOPT}.
In Figure~\ref{plotX} we show the results for the state function $x$
and in Figure~\ref{plotU} we show the results for the control function $u$,
where we have used $K=3$ and $n=100$.


\section*{Acknowledgments}

This article was supported by Portuguese funds through the
\emph{Center for Research and Development in Mathematics and Applications} (CIDMA),
and \emph{The Portuguese Foundation for Science and Technology} (FCT),
within project PEst-OE/MAT/UI4106/2014. Torres was also supported
by the FCT project PTDC/EEI-AUT/1450/2012,
co-financed by FEDER under POFC-QREN
with COMPETE reference FCOMP-01-0124-FEDER-028894.



\end{document}